\documentclass[11pt,oneside]{article}

\usepackage{amsfonts}
\usepackage[mathscr]{euscript}

\newcommand{\lbl}[1]{\label{#1}}

\newtheorem{theo}{Theorem}[section]
\newtheorem{prop}{Proposition}[section]
\newtheorem{lem}{Lemma}[section]
\newtheorem{remark}{Remark}[section]

\newcommand{\be}{\begin{equation}}
\newcommand{\ee}{\end{equation}}
\newcommand\bes{\begin{eqnarray}} \newcommand\ees{\end{eqnarray}}
\newcommand\bedd{\bes\left\{\begin{array}{ll}\medskip}
\newcommand\eedd{\end{array}\right.\ees}
\newcommand{\bess}{\begin{eqnarray*}}
\newcommand{\eess}{\end{eqnarray*}}

\newcommand\ep{\varepsilon}
\newcommand\lk{\left}

\newcommand\bbb{\big}
\newcommand\ff{, \ \ \forall \ }
\newcommand\nm{\nonumber}
\newcommand\dd{\displaystyle}
\newcommand\vp{\varphi}

\newcommand\ud{\underline}

\begin{document}
  \pagestyle{myheadings}

\begin{center}{\Large\bf On some free boundary problems of the prey-predator model}\footnote{This work was
supported by NSFC Grants 11071049 and 11371113}\\[4mm]
 {\Large  Mingxin Wang\footnote{{\sl E-mail}: mxwang@hit.edu.cn}}\\[1mm]
{\small  Natural Science Research Center, Harbin Institute of Technology, Harbin 150080, PR China}
\end{center}

\begin{quote}
\noindent{\bf Abstract.} In this paper we investigate some free boundary problems for the Lotka-Volterra type prey-predator model in one space dimension. The main objective is to understand the asymptotic behavior of the two species (prey and predator) spreading via a free boundary. We prove a spreading-vanishing dichotomy, namely the two species either successfully spread to the entire space as time $t$ goes to infinity and survive in the new environment, or they fail to establish and die out in the long run. The long time behavior of solution and criteria for spreading and vanishing are also obtained. Finally, when spreading successfully, we provide an estimate to show that the spreading speed (if exists) cannot be faster than the minimal speed of traveling wavefront solutions for the prey-predator model on the whole real line without a free boundary.

\noindent{\bf Keywords:} Prey-predator model; Free boundary problems; Spreading-vanishing dichotomy; Long time behavior; Criteria.

\noindent {\bf AMS subject classifications (2000)}:
35K51, 35R35, 92B05, 35B40.
 \end{quote}

 \section{Introduction}
 \setcounter{equation}{0} {\setlength\arraycolsep{3pt}

Understanding of spatial and temporal behaviors of interacting species in ecological systems is a central issue in population ecology. One aspect
of great interest for a model with multispecies interactions is whether the species can spread successfully. A lot of mathematicians have made
efforts to develop various models and investigated them from a viewpoint
of mathematical ecology. In this paper we consider three free boundary
problems for the Lotka-Volterra type prey-predator model.

In the real world, the following phenomenona will happen constantly:

(i)\, At the initial state, one kind of pest species (prey) occupied some bounded area (initial habitat). In order to control the pest species, the biological method is to put some natural enemies (predator) in this area;

(ii)\, There is some kind of species (prey) in a bounded area (initial habitat), and at some time (initial time) another kind of species (the new or invasive species, predator) enters this area.

In general, both prey and predator have a tendency to emigrate from the boundary to obtain their new habitat, i.e., they will move outward along the unknown curve (free boundary) as time increases. We argue that such prey in this model is the most favored food of the predator as causes of the {\it features of partial eclipse and picky eaters} for species, and its spreading behavior has such a dominant influence of spreading
of the predator that they roughly share the same spreading front. It is reasonable to assume that the free boundary invades at a rate that is proportional to the magnitudes of the prey and predator populations gradients there. We want to realize the dynamics/variations of  prey, predator and free boundary. For simplicity, we assume that the interaction between prey and predator obeys the Lokta-Volterra law, and focus on our problem to the one dimensional case. Under the suitable rescaling, the model we are concerned here becomes one of the following three free boundary problems:
\vspace{-2mm}\begin{quote}
(i) Left boundary is fixed with zero Dirichlet boundary conditions and the right boundary is free
 $$
 \left\{\begin{array}{lll}
 u_t-u_{xx}=u(1-u-av), &t>0,\ \ 0<x<h(t),\\[1mm]
  v_t-dv_{xx}=v(c-v+bu),\ \ &t>0, \ \ 0<x<h(t),\\[1mm]
 u=v=0,\ \ \ &t>0, \ \ x=0,\\[1mm]
 u=v=0,\ \ h'(t)=-\mu\big(u_x+\rho v_x\big), \ &t>0,\ \ x=h(t),\\[1mm]
 u(0,x)=u_0(x), \ \ v(0,x)=v_0(x),&x\in [0,h_0],\\[1mm]
 h(0)=h_0;
 \end{array}\right.\eqno{\rm(DFB)} $$
 (ii) Left boundary is fixed with zero Neumann boundary conditions and the right boundary is free
 $$
 \left\{\begin{array}{lll}
 u_t-u_{xx}=u(1-u-av), &t>0,\ \ 0<x<h(t),\\[1mm]
  v_t-dv_{xx}=v(c-v+bu),\ \ &t>0, \ \ 0<x<h(t),\\[1mm]
 u_x=v_x=0,\ \ \ &t>0, \ \ x=0,\\[1mm]
 u=v=0,\ \ h'(t)=-\mu\big(u_x+\rho v_x\big), \ &t>0,\ \ x=h(t),\\[1mm]
  u(0,x)=u_0(x), \ \ v(0,x)=v_0(x),&x\in [0,h_0],\\[1mm]
 h(0)=h_0;
 \end{array}\right.\eqno{\rm(NFB)} $$
(iii) With double free boundaries
 $$ \left\{\begin{array}{lll}
 u_t-u_{xx}=u(1-u-av), &t>0,\ \ g(t)<x<h(t),\\[1mm]
  v_t-dv_{xx}=v(c-v+bu),\ \ &t>0, \ \ g(t)<x<h(t),\\[1mm]
 u=v=0,\ \ g'(t)=-\mu_l\big(u_x+\rho_lv_x\big), \ &t>0, \ \ x=g(t),\\[1mm]
 u=v=0,\ \ h'(t)=-\mu_r\big(u_x+\rho_rv_x\big), \ &t>0,\ \ x=h(t),\\[1mm]
 u(0,x)=u_0(x), \ \ v(0,x)=v_0(x),&x\in [-h_0,h_0],\\[1mm]
 g(0)=-h_0,\ \ h(0)=h_0.
 \end{array}\right.\eqno{\rm(TFB)}$$
 \end{quote}\vspace{-2mm}
In the above three problems, $x=g(t)$ and $x=h(t)$ represent the left and
right moving boundaries, respectively, which are to be determined, $a, b, c, d, h_0, \mu, \rho, \mu_l, \mu_r, \rho_l$ and $\rho_r$ are given positive constants. The initial functions $u_0(x),v_0(x)$ satisfy

(DFB1) $u_0,\,v_0\in C^2([0,h_0])$, $u_0(0)=v_0(0)=u_0(h_0)=v_0(h_0)=0$, $u_0(x), v_0(x)>0$ in  $(0,h_0)$ for the problem (DFB);

(NFB1) $u_0,\,v_0\in C^2([0,h_0])$, $u_0'(0)=v_0'(0)=u_0(h_0)=v_0(h_0)=0$, $u_0(x), v_0(x)>0$ in  $(0,h_0)$ for the problem (NFB);

(TFB1) $u_0,\,v_0\in C^2([-h_0,h_0])$, $u_0(\pm h_0)=v_0(\pm h_0)=0$, $u_0(x), v_0(x)>0$ in  $(-h_0,h_0)$ for the problem (TFB).

In both problems (DFB) and (NFB), it is assumed that the species
can only invade further into the environment from the right end of the initial region. While in the proble (TFB),  it is assumed that the species can invade further into the environment from two ends of the initial region.

The ecological backgrounds of free boundary conditions in the above problems
can also refer to \cite{BDK}. Such kind of free boundary conditions has been used in \cite{GW, HMS}, \cite{Lin}--\cite{MYS3} and \cite{WZ}.

Recently, Wang and Zhao \cite{WZ, ZW} studied the similar free boundary problems to (TFB) with double free boundaries in which the prey lives in the whole space but the predator lives in the region enclosed by the free bounadry. Especially, in \cite{ZW}, the authors dealt with the higher dimension and heterogeneous environment case.  They have established the spreading-vanishing dichotomy, long time behavior of solution and criteria for spreading and vanishing.

Since the solution of (TFB) has the same properties as that of (NFB), we only discuss problems (DFB) and (NFB) in the following. For the global existence, uniqueness and estimates of solution, in the same way as \cite{ChenA, DLin} (see also \cite[Theorem 2.1, Lemma 2.1 and Theorem 2.2]{WZ}) we can prove the following theorem.

\begin{theo}\lbl{th1.1} \ Any one of {\rm(DFB)} and {\rm(NFB)} has a unique global solution, and for any $\alpha\in(0,1)$ and $T>0$,
  $$(u,v,h)\in \big[C^{\frac{1+\alpha}2,1+\alpha}(\overline{D}_T)\big]^2\times
 C^{1+\frac \alpha 2}([0,T]),$$
where
 \[D_T=\lk\{(t,x)\in \mathbb{R}^2:\,t\in(0,T],\, x\in\big(0,\,h(t)\big)\right\}.\]
Furthermore, there exists a positive constant $M$ such that
  \bess
&0<u(t,x), v(x,t)\leq M\ \ {\rm for}\ 0< t<\infty,\ 0<x<h(t),&\\[1mm]
 &0<h'(t)\leq M\ \ \, {\rm for}\ 0<t<\infty.&\eess
\end{theo}

In the absence of $v$, problems (DFB) and (NFB) are reduced to the one phase Stefan problems which were studied by Kaneko and Yamada \cite{KY} and Du and Lin \cite{DLin}, respectively. Free boundary problems for the logistic type model, including the higher dimension case, heterogeneous environment case (variable coefficients, time-periodic environment) and with seasonal succession, has been studied by many authors; please refer to, for example, \cite{DG, DG1, DGP, DLiang, DLou, DMZ, PZ}. There is a vast literature on the Stefan problems, and some important theoretical advances can be found in monographs \cite{CS, CR, RU} and the references therein. Some similar free boundary problems have been used in two-species models in several earlier papers; please refer to, for example, \cite{HIMN, HMS, Lin, MYS1, MYS2, MYS3} over a bounded spatial interval, and \cite{GW} over the half spatial line for the competition model.

The organization of this paper is as follows. To study the long time behavior of solution to the problem (DFB), in Section $2$ we discuss its stationary solutions. Section $3$ is devoted to the long time behavior of $(u,v)$ and get a spreading-vanishing dichotomy. To establish the criteria for spreading and vanishing, in Section $4$ we provide some comparison principles. The criteria for spreading and vanishing will be given in Section $5$. In Section $6$, we study the estimation of asymptotic spreading speed. The last section is a brief discussion.

\section{Positive solutions of the corresponding elliptic problems in half line}
 \setcounter{equation}{0} {\setlength\arraycolsep{3pt}

To discuss the long time behaviour of solution to the problem (DFB), we first discuss its stationary solutions. The stationary problem of (DFB) is the following elliptic problem in the half line
\bes\left\{\begin{array}{ll}\medskip
 -u''=u(1-u-av), &0<x<\infty,\\\medskip
 -dv''=v(c-v+bu),\ \ &0<x<\infty, \\
 u(0)=v(0)=0.
 \lbl{2.1}\end{array}\right.\ees
The main purpose of this section is to study the existence of positive solutions of (\ref{2.1}). To this aim, we first study the existence and uniqueness of positive solution to the following problem of single equation:
 \bes\left\{\begin{array}{ll}\medskip
 -du''=u\big(f(x)-\lambda u\big), \ \ 0<x<\infty,\\
 u(0)=0,
 \lbl{2.2}\end{array}\right.\ees
where $d$ and $\lambda$ are positive constants. When $f$ is a positive constant,
by Proposition 2.1 of \cite{BDK} or Proposition 4.1 of \cite{DLin}, problem (\ref{2.2})
has a unique positive solution $u(x)$. Moreover, $u'(x)>0$ and $\dd\lim_{x\to\infty}u(x)=f/\lambda$.

\subsection{The existence, uniqueness and stability of positive solution to (\ref{2.2})}

\begin{theo}\lbl{th2.1} \ Assume that $f$ satisfies
 \bes
 f\in C^\alpha_{\rm loc}([0,\infty)) \ \ \mbox{with} \ 0<\alpha<1, \ \
 \inf_{x\geq 0}f(x):=f_0>0, \ \ \|f\|_\infty<\infty.
 \lbl{2.3}\ees
Then the problem $(\ref{2.2})$ has a unique positive solution $u(x)$. Furthermore,

{\rm(i)}\, if $f(x)$ is increasing in $x$, so is $u(x)$ and $\lim\limits_{x\to\infty}u(x)=
\frac 1\lambda\lim\limits_{x\to\infty}f(x)$;

\vskip 4pt {\rm(ii)}\, if $f(x)$ is decreasing in $x$, then either $u(x)$ is increasing in $x$, or there exists $x_0>0$ such that $u(x)$ is increasing in $(0,x_0)$ and $u(x)$ is decreasing in $(x_0,\infty)$. Therefore, $\lim\limits_{x\to\infty}u(x)=\frac 1\lambda\lim\limits_{x\to\infty}f(x)$.
\end{theo}

{\bf Proof}.  We first analyse the properties of positive solution $u$ to the problem (\ref{2.2}). It is obvious that $u'(0)>0$. Moreover,
 \vspace{-2mm}\begin{quote}
(a)\, if $u(x_1)>\frac 1\lambda\|f\|_\infty$ and $u'(x_1)\geq 0$ for some $x_1\in(0,\infty)$, then $u(x)$ must approach infinity as $x$ approaches some finite $x_2$. This is impossible since $u(x)$ is defined in $[0,\infty)$; \\
(b)\, if $u(x_1)<f_0/\lambda$ and $u'(x_1)\leq 0$ for some $x_1\in(0,\infty)$, then $u$ must vanish at some finite $x_2$. This is also impossible as $u(x)$ is positive in $(0,\infty)$.
 \end{quote}
So we have
 \bes
 u'(0)>0, \ \ \ \mbox{and} \ \ \sup_{x\geq 0}u(x)\leq\frac 1\lambda\|f\|_\infty,
 \lbl{2.4}\ees
and there exists some positive constant $\tau$, depending on $d$, $\lambda$ and $f$, such that
 \bes
 u(x)\geq\tau , \ \ \ \forall \ x\geq 1.
 \lbl{2.5}\ees

Next, we prove the existence. It is well known that, when $l>\pi\sqrt{d/f_0}$, the problem
  \bes\left\{\begin{array}{ll}\medskip
 -du''=u\big(f(x)-\lambda u\big), \ \ 0<x<l,\\
 u(0)=0=u(l)
 \lbl{2.6}\end{array}\right.\ees
has a uniqueness positive solution, denoted by $u_l$, and $u_l$ satisfies
 $$\sup_{0\leq x\leq l}u_l(x)\leq\frac 1\lambda\|f\|_\infty.$$
Moreover, by the comparison principle
(\cite[Lemma 5.2]{Du}, or \cite[Proposition 2.1]{LPW}) we have that
$u_l(x)$ is increasing in $l$. In view of the regularity theory and compactness argument,
it follows that there exists a positive function $u$, such that $u_l\longrightarrow u$
in $C^2_{\rm loc}([0,\infty))$ as $l\longrightarrow\infty$, and $u$ solves (\ref{2.2}).

The uniqueness is followed from the following comparison principle (Proposition \ref{p2.1}).

At last, notice that for any $x\in(0,\infty)$, $u(x)\leq\frac 1\lambda f(x)$ if $u''(x)\leq 0$, while  $u(x)\geq\frac 1\lambda f(x)$ if $u''(x)\geq 0$, it is easily seen that conclusions (i) and (ii) hold.  \ \ \ \fbox{}

\vskip 3pt Now we give the comparison principle. Let $d$ and $\lambda$ be positive constants. Assume that $f_i\in C^\alpha_{\rm loc}([0,\infty))$ and satisfies $0<\inf_{x\geq 0}f_i(x)\leq\|f_i\|_\infty<\infty$ for $i=1,2$. By Theorem \ref{th2.1}, the problem
 \bes\left\{\begin{array}{ll}\medskip
 -du''=u\big(f_i(x)-\lambda u\big), \ \ 0<x<\infty,\\
 u(0)=0
 \nm\end{array}\right.\ees
has a positive solution, denoted by $u_i$.

\begin{prop}{\rm(}Comparison principle{\rm)}\lbl{p2.1} \ Under the above conditions,
if $f_1(x)\leq f_2(x)$ for all $x\geq 0$, then we have that
 \[u_1(x)\leq u_2(x), \ \ \ \forall \ x\geq 0.\]
 \end{prop}

{\bf Proof}. First, by (\ref{2.4}) and (\ref{2.5}), we have that $u_i'(0)>0$ and
$0<u_i(x)\leq\frac 1\lambda \sup_{x\geq 0}f_i(x)<\infty$ for all $x>0$, and $u_i(x)\geq\tau$
for all $x\geq 1$ and some positive constant $\tau$. Hence, there exists a constant $k\geq 1$
such $u_1(x)\leq ku_2(x)$ for all $x\geq 0$. Let
 \[k_0=\inf\big\{k>0: u_1(x)\leq ku_2(x), \ \ \forall \ x\geq 0\big\}, \ \ \mbox{i.e.}, \ \
 k_0=\sup_{x\geq 0}\frac{u_1(x)}{u_2(x)}.\]
Then $u_1(x)\leq k_0u_2(x)$ for all $x\geq 0$. If $k_0\leq 1$, then our conclusion is true.

We assume that $k_0>1$ and shall derive a contradiction. Let $\vp(x)=k_0u_2(x)-u_1(x)$.
Then $\vp(x)\geq 0$ for all $x\geq 0$, and $\vp$ satisfies
 \bes
  -d\vp''&=&k_0u_2\big(f_2(x)-\lambda u_2\big)-u_1(x)\big(f_1(x)-\lambda u_1\big)\nm\\
  &\geq&f_1(x)\vp-\lambda(k_0u_2^2-u_1^2)\nm\\
  &=&\big[f_1(x)-\lambda(k_0u_2+u_1)\big]\vp+\lambda k_0u_2^2(k_0-1), \ \ 0<x<\infty \lbl{2.7}\\
  &>&\big[f_1(x)-\lambda(k_0u_2+u_1)\big]\vp\nm\\
  &\geq&\frac 12 f_1(0)\vp, \ \ 0<x\ll 1\lbl{2.8}\ees
since $f_1(0)>0$ and $u_i(x)\longrightarrow 0$ as $x\longrightarrow 0$. Thanks to $\vp(0)=0$, we have $\vp'(0)\geq 0$. If $\vp'(0)=0$, it follows from (\ref{2.8}) that $\vp'(x)<0$ as $0<x\ll 1$. Consequently, $\vp(x)<0$ as $0<x\ll 1$. This is a contradiction.
So, $\vp'(0)>0$, i.e., $k_0u_2'(0)>u_1'(0)$. Remember this fact and the definition of $k_0$, it is easily seen that at least one of the following happens:

  (i)\ There exists some $x_0\in(0,\infty)$ such that $k_0u_2(x_0)=u_1(x_0)$,

 (ii)\, $k_0u_2(x)-u_1(x)>0$ for all $x>0$ and $\dd\liminf_{x\to\infty}\big[k_0u_2(x)-u_1(x)\big]=0$.

\vskip 4pt When the case (i) occurs, then $\vp(x_0)=0$, $\vp'(x_0)=0$ and $\vp''(x_0)\geq 0$. It is derived from (\ref{2.7}) that $k_0\leq 1$. This is a contradiction.

When the case (ii) happens, then $\dd\liminf_{x\to\infty}\vp(x)=0$. Remember that $\vp(x)>0$ for all $x>0$, it is not hard to prove that there exists a sequence $\{x_n\}$ with $x_n\longrightarrow\infty$ such that
 \[\vp''(x_n)\geq 0, \ \ \ \lim_{n\to\infty}\vp(x_n)=0.\]
By passing to a subsequence, we may assume that $u_2(x_n)\longrightarrow \sigma$ for some positive constant $\sigma$. It follows from (\ref{2.7}) that
 \[0\geq\big[f_1(x_n)-\lambda\bbb(k_0u_2(x_n)+u_1(x_n)\big)\big]\vp(x_n)
 +\lambda k_0u_2^2(x_n)(k_0-1).\]
Letting $n\longrightarrow\infty$ we get
 \[0\geq \lambda k_0\sigma^2(k_0-1)>0.\]
This is also a contradiction. The proof is finished.  \ \ \ \fbox{}

In the following, we denote the unique positive solution of (\ref{2.2}) by $\hat u(x)$.
Now we discuss its global stability.

\begin{theo}\lbl{th2.2} \ Assume that $d$ and $\lambda$ are positive constants, and $f(x)$ satisfies {\rm(\ref{2.3})}. Suppose that $\phi\not\equiv 0$ is a bounded, continuous and nonnegative function. Let $u(t,x)$ be the unique solution of the following parabolic problem
 \bes\left\{\begin{array}{ll}\medskip
 u_t-du_{xx}=u\big(f(x)-\lambda u\big), \ \ &t>0, \ 0<x<\infty,\\\medskip
 u(t,0)=0, &t>0,\\
 u(0,x)=\phi(x), &0\leq x<\infty.
 \nm\end{array}\right.\ees
Then
 \bes
 \lim_{t\to\infty}u(t,x)=\hat u(x)\ \ \ \mbox{uniformly in any compact
 subset of } \, [0, \infty).
 \lbl{2.9}\ees
 \end{theo}

{\bf Proof}. \ By the positivity of parabolic equations we have that $u(t,x)>0$ for all $t>0$ and $x>0$, and $u_x(t,0)>0$ for all $t>0$. We may assume that $\phi(x)>0$ for all $x>0$, and $\phi(0)=0$ and $\phi'(0)>0$. It is well known that, when $l>\pi\sqrt{d/f_0}$, where $f_0=\inf_{x\geq 0}f(x)>0$, the problem (\ref{2.6}) has a unique positive solution, denoted by $\hat u_l(x)$. For $l>\pi\sqrt{d/f_0}$, consider the following initial-boundary value problem
   \bes\left\{\begin{array}{ll}\medskip
 u_t-du_{xx}=u\big(f(x)-\lambda u\big), \ \ &t>0, \ 0<x<l,\\\medskip
 u(t,0)=0=u(t,l), &t>0,\\
 u(0,x)=\phi(x), &0\leq x\leq l.
 \lbl{2.10}\end{array}\right.\ees
Let $u_l(t,x)$ be the unique solution of (\ref{2.10}). Then
 \bes
 u(t,x)\geq u_l(t,x), \ \ \ \forall \ x\in[0,l], \ t\geq 0.
 \lbl{2.11}\ees
It is obvious that a large positive constant $M$ is an upper solution of (\ref{2.6}). Note that $l>\pi\sqrt{d/f_0}$, it can be proved that $\delta_1\sin\frac{\pi x}l$ is a lower solution of (\ref{2.6}) provided that $0<\delta_1\ll 1$. Since $\phi(x)>0$ for $x>0$ and $\phi'(0)>0$, there is a small positive constant $\delta_2$ such that $\delta_2\sin\frac{\pi x}l\leq\phi(x)$ for all $x\in[0,l]$. Set $\delta=\min\{\delta_1,\delta_2\}$, then  $\delta\sin\frac{\pi x}l\leq\phi(x)$ for all $x\in[0,l]$ and $\delta\sin\frac{\pi x}l$ is a lower solution of (\ref{2.6}). Let $\bar u_l(t,x)$ and $\ud u_l(t,x)$ be the unique solution of (\ref{2.10}) with $\phi(x)=M$ and $\phi(x)=\delta\sin\frac{\pi x}l$, respectively. Then
 \bes
 \ud u_l(t,x)\leq u_l(t,x)\leq\bar u_l(t,x), \ \ \ \forall \ x\in[0,l], \ t\geq 0,
 \lbl{2.12}\ees
and $\bar u_l(t,x)$ is decreasing and $\ud u_l(t,x)$ is increasing in $t$. Moreover, both limits $\lim_{t\to\infty}\bar u_l(t,x)=\bar u_l(x)$ and $\lim_{t\to\infty}\ud u_l(t,x)=\ud u_l(x)$ are positive solutions of (\ref{2.6}). By the uniqueness we have $\bar u_l(x)\equiv\ud u_l(x)\equiv \hat u_l(x)$. It follows from (\ref{2.12}) that
 \bes
 \lim_{t\to\infty}u_l(t,x)=\hat u_l(x)\ \ \ \mbox{uniformly in } \, [0, l].
 \nm\ees
Hence, by (\ref{2.11}),
   \bes
 \liminf_{t\to\infty}u(t,x)\geq\hat u_l(x)\ \ \ \mbox{uniformly in } \, [0, l].
 \lbl{2.13}\ees
From the proof of Theorem \ref{th2.1} we know that $\hat u_l\longrightarrow\hat u$ in $C^{2+\alpha}_{\rm loc}([0,\infty))$ as $l\longrightarrow\infty$. Consequently, by (\ref{2.13}), we arrive at
 \bes\liminf_{t\to\infty}u(t,x)\geq\hat u(x)\ \ \ \mbox{uniformly in any compact
 subset of } \, [0, \infty).
 \lbl{2.14}\ees

Let $M=\max\lk\{\|\phi\|_\infty, \, \|f\|_\infty/\lambda\right\}$ and $\bar u(t,x)$ be the unique solution of \bes\left\{\begin{array}{ll}\medskip
 u_t-du_{xx}=u\big(f(x)-\lambda u\big), \ \ &t>0, \ 0<x<\infty,\\\medskip
 u(t,0)=0, &t>0,\\
 u(0,x)=M, &0\leq x<\infty.
 \nm\end{array}\right.\ees
Then we have the following conclusions:
  \vspace{-2mm}\begin{quote}
(a)\, $\bar u(t,x)$ is monotone decreasing in $t$;\\
(b)\, $u(t,x)\leq\bar u(t,x)$ for $x\geq 0$ and $t\geq 0$;\\
(c)\, $\bar u(t,x)\geq\hat u(x)$ for $x\geq 0$ and $t\geq 0$ since $\hat u(x)<M$;\\
(d)\, $\lim_{t\to\infty}\bar u(t,x)=u^*(x)$ uniformly in any compact subset of $[0, \infty)$ for some positive solution $u^*(x)$ of (\ref{2.2}).
 \vspace{-2mm} \end{quote}
Because $\hat u(x)$ is the unique positive solution of (\ref{2.2}), one has $u^*(x)=\hat u(x)$. Therefore,
\bes\limsup_{t\to\infty}u(t,x)\leq\hat u(x)\ \ \ \mbox{uniformly in any compact
 subset of } \, [0, \infty).
 \nm\ees
This combines with (\ref{2.14}) to derive (\ref{2.9}). \ \ \ \fbox{}

Before ending this subsection, we provide two propositions that will be used in the study of the long time behaviour of solution to the problem (DFB).

\begin{prop}\lbl{p2.2} \ Assume that $d$ and $\lambda$ are positive constants, $f(x)$ satisfies {\rm(\ref{2.3})}. For any given constant $K>\frac1\lambda\|f\|_\infty$ and any $l\gg 1$, let $\bar u_l(x)$ be the unique positive solution of
 \bes\left\{\begin{array}{ll}\medskip
 -du''=u\big(f(x)-\lambda u\big), \ \ 0<x<l,\\
 u(0)=0, \ \ \ u(l)=K.
 \nm\end{array}\right.\ees
Then we have
 \[\lim_{l\to\infty}\bar u_l(x)=\hat u(x)\ \ \ \mbox{uniformly in any compact
 subset of } \, [0, \infty).\]
 \end{prop}

{\bf Proof}. \ First, the existence and uniqueness of $\bar u_l(x)$ can be obtained by the upper and lower solutions method and the comparison principle (\cite[Lemma 5.2]{Du}, \cite[Proposition 2.1]{LPW}), respectively.

Since $\bar u_l(x)\leq K$ in $[0,l]$, by the comparison principle we see that $\bar u_l(x)\geq \hat u(x)$ in $[0,l]$, and $\bar u_l(x)$ is decreasing in $l$. Similar to the proof of Theorem \ref{th2.1}, $\bar u_l\longrightarrow\hat u$ in $C^2_{\rm loc}([0,\infty))$ as $l\longrightarrow\infty$ because $\hat u$ is the unique positive solution of (\ref{2.2}). \ \ \ \fbox{}

Using Proposition \ref{p2.1}, and the regularity theory and compactness argument, we can prove the following proposition.

\begin{prop}\lbl{p2.3} \ Assume that $d$ and $\lambda$ are positive constants, $f(x)$ satisfies {\rm(\ref{2.3})}. Let $0<\ep\ll 1$ and $u_\ep^{\pm}(x)$ be the unique positive solution of
 \bes\left\{\begin{array}{ll}\medskip
 -du''=u\big(f(x)\pm\ep-\lambda u\big), \ \ 0<x<\infty,\\
 u(0)=0.
 \nm\end{array}\right.\ees
Then
 \[\lim_{\ep\to 0}\bar u_\ep^{\pm}(x)=\hat u(x)\ \ \ \mbox{uniformly in any compact
 subset of } \, [0, \infty).\]
 \end{prop}

\subsection{The existence of positive solution to (\ref{2.1})}

\begin{theo}\lbl{th2.3} \ Assume that $a(b+c)<1$. Then the problem {\rm(\ref{2.1})}
has a positive solution. Moreover, any positive solution $(u,v)$ of {\rm(\ref{2.1})} satisfies
\bes
 \ud u(x)\leq u(x)\leq\bar u(x), \ \ \ \ud v(x)\leq v(x)\leq\bar v(x), \ \ \
 \forall \ x\geq 0,\lbl{2.15}\ees
where $\bar u,\,\bar v,\,\ud u$ and $\ud v$ will be given in the following proof.
\end{theo}

{\bf Proof}. {\it Step 1}: The construction of $\ud u$, $\ud v$, $\bar u$ and $\bar v$.

Let $\bar u$ be the unique positive solution of
 \bes\left\{\begin{array}{ll}\medskip
 -u''=u(1-u), \ \ 0<x<\infty,\\
 u(0)=0.
 \lbl{2.16}\end{array}\right.\ees
Then
 \[\bar u'(x)>0, \ \ \mbox{and} \ \ \lim_{x\to\infty}\bar u(x)=1.\]
In view of Theorem \ref{th2.1}, the problem
 \bes\left\{\begin{array}{ll}\medskip
 -dv''=v\big[c-v+b\bar u(x)\big],\ \ &0<x<\infty, \\
 v(0)=0
 \lbl{2.17}\end{array}\right.\ees
has a unique positive solution, denoted by $\bar v(x)$. Then $\bar v'(x)>0$ and $\dd\lim_{x\to\infty}\bar v(x)=b+c$. Since $a(b+c)<1$, by Theorem \ref{th2.1} again, the problem
 \bes\left\{\begin{array}{ll}\medskip
 -u''=u\big(1-a\bar v(x)-u\big), \ \ 0<x<\infty,\\
 u(0)=0
 \nm\end{array}\right.\ees
has a unique positive solution, denoted by $\ud u(x)$. Then $\ud u(x)<1$ and  $\dd\lim_{x\to\infty}\ud u(x)=1-a(b+c)$. By virtue of Theorem \ref{th2.1} once again, the problem
 \bes\left\{\begin{array}{ll}\medskip
 -dv''=v\big[c-v+b\ud u(x)\big],\ \ &0<x<\infty, \\
 v(0)=0
 \nm\end{array}\right.\ees
has a unique positive solution, denoted by $\ud v(x)$. Then $\ud v(x)<b+c$.

Applying Proposition \ref{p2.1} we have that $\ud u(x)\leq \bar u(x)$ for all $x\geq 0$. Consequently,  $\ud v(x)\leq \bar v(x)$ for all $x\geq 0$ by use of Proposition \ref{p2.1} once again.

{\it Step 2}: Existence of positive solution.

The conclusion of Step $1$ shows that $\ud u$, $\ud v$, $\bar u$ and $\bar v$ are the coupled ordered lower and upper solutions of (\ref{2.1}). For any given $l>0$, it is obvious that $\ud u$, $\ud v$, $\bar u$ and $\bar v$ are also the coupled ordered lower and upper solutions of the following problem
 \bess\left\{\begin{array}{ll}\medskip
 -u''=u(1-u-av), &0<x<l,\\\medskip
 -dv''=v(c-v+bu),\ \ &0<x<l, \\
 u(0)=\ud u(0), \ \ v(0)=\ud v(0), \\[2mm]
 u(l)=\bar u(l), \ \ v(l)=\bar v(l).
 \end{array}\right.\eess
By the standard upper and lower solutions method we know that the above problem has at  least one positive solution, denoted by $(u_l,v_l)$, and
 \[\ud u(x)\leq u_l(x)\leq\bar u(x), \ \ \ \ud v(x)\leq v_l(x)\leq\bar v(x), \ \ \
 \forall \ 0\leq x\leq l.\]
Applying the local estimation and compactness argument, it can be concluded that there exists a pair $(u,v)$, such that $(u_l,v_l)\longrightarrow (u,v)$ in $\lk[C^2_{\rm loc}([0,\infty))\right]^2$, and $(u, v)$ solves (\ref{2.1}). It is obvious that $(u,v)$ satisfies (\ref{2.15}).

By use of Proposition \ref{p2.1} successively, we also see that if $a(b+c)<1$ then any positive solution $(u,v)$ of (\ref{2.1}) satisfies (\ref{2.15}).  \ \ \ \fbox{}

\section{Long time behavior of $(u,v)$ and spreading-vanishing dichotomy}
\setcounter{equation}{0}{\setlength\arraycolsep{2pt}

The main purpose of this section is to study the long time behavior of $(u,v)$ and obtain a spreading-vanishing dichotomy. It follows from Theorem \ref{th1.1} that $x=h(t)$ is
monotonic increasing. Therefore, $\lim_{t\to\infty} h(t)=h_\infty\in(0,\infty]$.
To discuss the long time behavior of $(u,v)$, we first derive an estimate.

\begin{theo}\lbl{th3.1} \ Let $(u,v,h)$ be the solution of {\rm(DFB)} or {\rm(NFB)}. If
$h_\infty<\infty$, then there exists a constant $M>0$, such that
 \[ \|u(t,\cdot), v(t,\cdot)\|_{C^1[0,\,h(t)]}\leq M,\ \forall \ t>1.\]
Moveover,
 \bess
 \lim_{t\to\infty}h'(t)=0.
 \eess
\end{theo}

{\bf Proof}. \ The proof is similar to that of Theorem 4.1 of \cite{WZ}. We omit the details.

\subsection{Vanishing case ($h_\infty<\infty$)}

To discuss the asymptotic behaviors of $u$ and $v$, we first prove a general result. Let $d$, $\beta$ and $g_0$ be positive constants and $C\in\mathbb{R}$. Assume that $w_0\in C^2([0,g_0])$ and satisfies $w_0'(0)=0$ {\rm(}or $w_0(0)=0${\rm)}, $w_0(g_0)=0$ and $w_0(x)>0$ in $(0,h_0)$. Let $$g\in C^{1+\frac \alpha 2}([0,\infty)), \ \ \ w\in C^{\frac{1+\alpha}2,1+\alpha}([0,\infty)\times[0,g(t)])$$
for some $\alpha>0$, and satisfy $g(t)>0, w(t,x)>0$ for all $0\leq t<\infty$ and $0<x<g(t)$.

\begin{prop}\label{p3.1} Under the above conditions, we further suppose that
  \bes
 \lim_{t\to\infty} g(t)=g_\infty<\infty, \ \ \lim_{t\to\infty} g'(t)=0,
 \label{3.1}\ees
and
 \be
 \|w(t,\cdot)\|_{C^1[0,\,g(t)]}\leq M, \ \ \forall \ t>1\label{3.2}
 \ee
for some constant $M>0$.  If $(w,g)$ satisfies
  \bess\left\{\begin{array}{lll}
 w_t-dw_{xx}\geq w(C-w), &t>0,\ \ 0<x<g(t),\\[1.5mm]
 w_x=0 \ (\mbox{or} \ w=0),\ \ \ &t>0, \ \ x=0,\\[1.5mm]
 w=0,\ \ g'(t)\geq-\beta w_x, \ &t>0,\ \ x=g(t),\\[1.5mm]
 w(0,x)=w_0(x), \ &x\in [0,g_0],\\[1.5mm]
 g(0)=g_0,
 \end{array}\right.\eess
then
  \be
 \lim_{t\to\infty}\,\max_{0\leq x\leq g(t)}w(t,x)=0.\label{3.3}
\ee
\end{prop}

{\bf Proof}. \  The proof is similar to that of \cite[Theorem 4.2]{WZ}. For the convenience to reader, we shall give the details because of this is a very important conclusion in the study of free boundary problems for systems.

On the contrary we assume that there exist $\varepsilon>0$ and $\left\{(t_j,x_j)\right\}_{j=1}^{\infty}$, with $0\leq x_j<g(t_j)$ and  $t_j\to\infty$ as $j\to\infty$, such that
 \bes w(t_j,x_j)\geq 3\varepsilon,\ j=1,2,\cdots.\label{3.4}\ees
Since $0\leq x_j<g_\infty$, there exist a subsequence of $\{x_j\}$, denoted by itself,
and $x_0\in[0,g_\infty]$, such that $x_j\to x_0$ as $j\to\infty$.
We affirm that $x_0<g_\infty$. If this is not true, then $x_j-g(t_j)\to 0$
as $j\to\infty$. By use of the inequality (\ref{3.4}) firstly and the inequality
(\ref{3.2}) secondly, we have that
 \bess
 \left|\dd\frac{4\varepsilon}{x_j-g(t_j)}\right|\leq\left|\frac{w(t_j,x_j)}{x_j-g(t_j)}\right|=
 \left|\frac{w(t_j,x_j)-w(t_j,g(t_j))}{x_j-g(t_j)}\right|=\left|w_x(t_j,\bar x_j)\right|\leq M,\eess
where $\bar x_j\in(x_j,g(t_j))$. It is a contradiction since $x_j-g(t_j)\to 0$. Similarly, it is easily deduced that $x_0\neq 0$ if the boundary condition at $x=0$ is $w=0$.

{\bf Case 1: $x_0>0$}.\, By use of (\ref{3.2}) and (\ref{3.4}), there exists $\delta:0<\delta<x_0$ such that $x_0+\delta<g_\infty$ and
 $$w(t_j,x)\geq 2\varepsilon,\ \ \ \forall\ x\in[x_0-\delta, \, x_0+\delta]$$
for all large $j$. Since $g(t_j)\to s_\infty$ as $j\to+\infty$, we may think that
$g(t_j)>x_0+\delta$ for all $j$. Let
 \[r_j(t)=x_0+\delta+t-t_j, \ \ \ \tau_j=\inf\left\{t>t_j:\, g(t)=r_j(t)\right\}.\]
Clearly, $x_0+\delta+\tau_j-t_j=r_j(\tau_j)<g_\infty$, $g(\tau_j)=r_j(\tau_j)$, and $x_0+\delta\leq r_j(t)\leq g(t)$ in $[t_j,\tau_j]$.
Define
 \bess
 y_j(t,x)&=&\dd\frac{2(x-x_0)-t+t_j}{2\delta+t-t_j}(\pi-\theta),\\[2mm]
  D_j&=&\{(t,x):t_j<t<\tau_j, \ x_0-\delta<x<r_j(t)\},
 \eess
and
  $$ w_j(t,x)=\varepsilon {\rm e}^{-k(t-t_j)}[\cos y_j(t,x)+\cos
\theta],\ \ \ (t,x)\in\overline D_j,$$
where $\theta\ (\theta<\pi/8)$ and $k$ are positive constants to be chosen
later. Then $|y_j(t,x)|\leq \pi-\theta$ for $(t,x)\in\overline D_j$, and
$y_j(t,x_0-\delta)=-(\pi-\theta)$, $y_j(t,r_j(t))=\pi-\theta$.
Therefore, $ w_j(t,x)\geq 0$ in $ D_j$, and  $ w_j(t,x_0-\delta)= w_j(t,r_j(t))=0$.

We want to compare $w(t,x)$ and $ w_j(t,x)$ in $D_j$. It is obvious that
 $$w(t,x_0-\delta)\geq 0= w_j(t,x_0-\delta), \ \ w(t,r_j(t))\geq 0=
 w_j(t,r_j(t)), \ \ \forall \  t\in[t_j,\tau_j],$$
and
 $$w(t_j,x)\geq 2\varepsilon> w_j(t_j,x), \ \ \forall \  x\in[x_0-\delta,\,x_0+\delta].$$

Thus, if the positive constants $\theta$ and $k$ can be chosen independent of
$j$ such that
 \bes w_{jt}-d w_{jxx}\leq w_j(C- w_j)
 \ \ \ \mbox{in} \ \ D_j,
 \label{3.7}\ees
then $ w_j(t,x)\leq w(t,x)$ in $ D_j$ by
applying the comparison principle to $w$ and $ w_j$. Once this is done, since $w(\tau_j,g(\tau_j))=0= w_j(\tau_j,r_j(\tau_j))$ and
$g(\tau_j)=r_j(\tau_j)$, we have $w_x(\tau_j,g(\tau_j))\leq
w_{jx}(\tau_j,r_j(\tau_j))$. Notice $\theta<\pi/8$ and $2\delta+\tau_j-t_j<g_\infty$, it follows that
  $$ w_{jx}(\tau_j,r_j(\tau_j))=-\dd\frac{2\varepsilon(\pi-\theta)}
{2\delta+\tau_j-t_j}{\rm e}^{-k(\tau_j-t_j)}\sin(\pi-\theta)
\leq-\dd\frac{7\varepsilon\pi}{4g_\infty}{\rm e}^{-kg_\infty}
 \sin\theta.$$
This fact combined with the boundary condition $-\beta w_x(\tau_j,g(\tau_j))\leq g'(\tau_j)$ provides us
  $$g'(\tau_j)\geq\dd\frac{7\beta\varepsilon\pi}{4g_\infty}
  {\rm e}^{-kg_\infty}\sin \theta,$$
which implies $\limsup_{t\to\infty}g'(t)>0$ as $\lim_{j\to\infty}\tau_j=\infty$. This contradicts to (\ref{3.1}), and (\ref{3.3}) is obtained.

In the following we shall demonstrate that (\ref{3.7}) holds as long as
$\theta$ and $k$ satisfy
\bes
 &\theta<\dd\frac{\pi}{8}, \ \ \ \sin\theta\leq
\dd\frac{3d\delta^2\pi}{g_\infty^3},&\label{3.8}\\[1.5mm]
 &k>2\varepsilon+|C|+d\left(\dd\frac{\pi}{\delta}\right)^2+\dd\frac{g_\infty\pi}{2\delta^2(\cos\theta-\cos
2\theta)}.&\label{3.9}
  \ees
To achieve this objective, remember $0\leq  w_j\leq 2\varepsilon$ and $2\delta+\tau_j-t_j<g_\infty$, a series of computations indicate that, for $(t,x)\in D_j$,
 \bess
&& w_{jt}-d w_{jxx}- w_j(C- w_j)\\[2mm]
&=&-k w_j-\varepsilon {\rm e}^{-k(t-t_j)}y_{jt}\sin
y_j+d\varepsilon {\rm e}^{-k(t-t_j)}y_{jx}^2\cos y_j-
w_j(C- w_j)\\[2mm]
&\leq&\dd\left(2\varepsilon +|C|+d y_{jx}^2-k\right)
w_j-d\varepsilon {\rm e}^{-k(t-t_j)}y_{jx}^2\cos
\theta-\varepsilon {\rm e}^{-k(t-t_j)}y_{jt}\sin y_j\\[2mm]
&\leq&\dd\left(2\varepsilon+|C|+d\left(\frac{\pi}{\delta}\right)^2-k\right) w_j
+\varepsilon{\rm e}^{-k(t-t_j)}\dd\left[\dd\frac{\pi(x-x_0+\delta)}{2\delta^2}|\sin y_j|-d\left(\dd\frac{2(\pi-\theta)}
{g_\infty}\right)^2\cos\theta\right]\\[1.5mm]
&:=&I(t,x).
  \eess
Clearly, $2\varepsilon+|C|+d\left(\frac{\pi}{\delta}\right)^2-k<0$ by (\ref{3.9}).   
Since $|y_j(t,x)|\leq \pi-\theta$ in $\overline D_j$, we can decompose $ D_j= D_j^1\bigcup
 D_j^2$, where
\bess
   D_j^1&=&\left\{(t,x)\in D_j:\, t_j<t<\tau_j,\
 \pi-2\theta<|y_j(t,x)|<\pi-\theta\right\},\\[1mm]
  D_j^2&=&\left\{(t,x)\in D_j:\, t_j<t<\tau_j,
 \ |y_j(t,x)|\leq\pi-2\theta\right\}.\eess
It is obvious that $|\sin y_j(t,x)|\leq \sin 2\theta$ in $ D_j^1$,
and $\cos y_j(t,x)\geq -\cos 2\theta$ in $ D_j^2$. Because $\theta<\pi/8$,  and $ w_j(t,x)\geq 0$, $0\leq  x-x_0+\delta\leq g_\infty$, $|\sin y_j(t,x)|\leq 1$ in $ D_j$, in view of (\ref{3.8}) and (\ref{3.9}), we conclude that
   \bess
   I(t,x)&<&\varepsilon{\rm e}^{-k(t-t_j)}\dd\left[\dd\frac{\pi(x-x_0+\delta)}{2\delta^2}|\sin y_j|-d\left(\dd\frac{2(\pi-\theta)}
 {g_\infty}\right)^2\cos\theta\right]\\[2mm]
 &\leq&\varepsilon {\rm e}^{-k(t-t_j)}\left(\dd\frac{\pi g_\infty}{2\delta^2}\sin 2\theta-d\dd\frac{3\pi^2}{g_\infty^2}\cos \theta\right)<0 \ \ \  {\rm in} \ \ D_j^1,\eess
and
 \bess
 I(t,x)&<&\dd\left[2\varepsilon+|C|+d\left(\frac{\pi}{\delta}\right)^2-k\right] w_j
+\varepsilon{\rm e}^{-k(t-t_j)}\frac{\pi(x-x_0+\delta)}{2\delta^2}|\sin y_j|\\[2mm]
 &\leq& \varepsilon {\rm e}^{-k(t-t_j)}\left\{\left[2\varepsilon+|C|
 +d\left(\frac{\pi}{\delta}\right)^2-k\right](\cos\theta-\cos 2\theta)
 +\dd\frac{g_\infty\pi}{2\delta^2}\right\}<0 \ \ \  {\rm in} \ \ D_j^2.\eess
And so, (\ref{3.7}) holds.

{\bf Case 2: $x_0=0$}.\, In this case, the boundary condition at $x=0$ must be $w_x=0$. Similar to the above, there exists $0<\delta<g_\infty$ such that, for all large $j$,
there holds $w(t_j,x)\geq 2\varepsilon$ in $[0,\,\delta]$. Let $r_j(t)=\delta+t-t_j$, and set
  $$\tau_j=\inf\left\{t>t_j:\, g(t)=r_j(t)\right\}, \ \ D_j=\{(t,x):t_j<t<\tau_j, \ 0<x<r_j(t)\}.$$
Define $y_j(t,x)=(\pi-\theta)\frac{x}{\delta+t-t_j}$,
 $$ w_j(t,x)=\varepsilon {\rm e}^{-k(t-t_j)}[\cos y_j(t,x)+\cos
\theta],\ \ \ (t,x)\in\overline D_j,$$
where $\theta$ and $k$ satisfy
\bess
 0<\theta<\dd\frac{\pi}{8}, \ \ \sin\theta\leq
\dd\frac{3\pi d\delta^2}{8g_\infty^3}, \ \
 k>2\varepsilon+|C|+d\left(\dd\frac{\pi}{\delta}\right)^2+\dd\frac{\pi g_\infty}{\delta^2(\cos\theta-\cos 2\theta)}.
  \eess
It is clear that $w_{jx}(t,0)=0=w_x(t,0)$, $w_j(t,r_j(t))=0\leq w(t,r_j(t))$ in  $[t_j,\tau_j]$, and $w_j(t_j,x)<2\varepsilon\leq w(t_j,x)$ in $[0,\,\delta]$.
The rest of the proof is similar to the case 1.\ \ \ \fbox{} 

\begin{theo}\lbl{th3.2} \ Let $(u,v,h)$ be any solution of {\rm(DFB)} or {\rm(NFB)}. If
$h_\infty<\infty$, then
 \bes
 \lim_{t\to\infty}\|u(t,\cdot),\,v(t,\cdot)\|_{C([0,\,h(t)])}=0.\lbl{3.10}\ees
This result shows that if both prey and predator can not spread into the infinity, then they
will die out eventually.
\end{theo}

We should remark that this theorem plays key roles in the following two aspects: (i) affirming that the two species disappear eventually; (ii) determining the criteria for spreading and vanishing (see the following Section $5$).

{\bf Proof of Theorem \ref{th3.2}}. \, Since $u(t,x)>0$ and $u_x(t,h(t))<0$, we see that $v$ satisfies
$$\left\{\begin{array}{lll}
 v_t-dv_{xx}\geq v(c-v),\ \ &t>0, \ \ 0<x<h(t),\\[1mm]
 v=0 \ (\mbox{or} \ v_x=0),\ \ \ &t>0, \ \ x=0,\\[1mm]
 v=0,\ \ h'(t)\geq-\mu\rho v_x, \ &t>0,\ \ x=h(t),\\[1mm]
 v(0,x)=v_0(x),&x\in [0,h_0],\\[1mm]
 h(0)=h_0.
 \end{array}\right.$$
In view of Theorem \ref{th3.1} and Proposition \ref{p3.1} we have that $\lim_{t\to\infty}\|v(t,\cdot)\|_{C([0,\,h(t)])}=0$. Hence, there exists $T\gg 1$, such that
 \[v(t,x)<1/(2a) \ \ \ \forall \ t\geq T, \ 0\leq x\leq h(t).\]
Remember that $v_x(t,h(t))<0$ for all $t\geq T$, we see that $u(t,x)$ satisfies
 \bess
 \left\{\begin{array}{lll}
  u_t-u_{xx}>u(1/2-u),\ \ &t\geq T, \ \ 0<x<h(t),\\[2mm]
 u=0 \ \ \mbox{or} \ \ u_x=0,\ \ &t\geq T, \ \ x=0,\\[2mm]
 u=0,\ \ h'(t)>-\mu u_x, \ &t\geq T,\ \ x=h(t),\\[2mm]
 u(t,x)=u(T,x),& t=T, \ \ x\in [0,\,h(t)].
 \end{array}\right.\eess
Applying Theorem \ref{th3.1} and Proposition \ref{p3.1} again, we conclude that
$\lim_{t\to\infty}\|u(t,\cdot)\|_{C([0,\,h(t)])}=0$. The proof is finished. \ \ \ \fbox{}

\subsection{Spreading case ($h_\infty=\infty$)}

We first consider the problem (NFB). In the same way as the proofs of Theorems 4.3 and 4.4 in \cite{WZ}, we can prove the following two theorems.

\begin{theo}\lbl{th3.3} Let $(u,v,h)$ be the solution of {\rm(NFB)}. If $h_\infty=\infty$, then for the weakly hunting case $ac<1$, $ab<1$, we have
 \bess
 \lim_{t\to\infty}u(t,x)=\frac{1-ac}{1+ab},\ \ \ \
\lim_{t\to\infty}v(t,x)=\frac{b+c}{1+ab}
  \eess
uniformly in any compact subset of $[0,\infty)$.
 \end{theo}

We remark that the conditions $ac<1$ and $ab<1$ are similar to the weak competition conditions, see \cite{GW}.

\begin{theo}\lbl{th3.4} Let $(u,v,h)$ be the solution of {\rm(NFB)}. If $h_\infty=\infty$, then for the strongly hunting case $ac\geq 1$, we have
 \[\lim_{t\to\infty}u(t,x)=0,\ \ \ \ \lim_{t\to\infty}v(t,x)=c\]
uniformly in any compact subset of $[0,\infty)$.
 \end{theo}

In the rest of this section, we consider the problem (DFB). This is the main part of this  section.

\begin{theo}\lbl{th3.5} \ Assume that $h(\infty)=\infty$. If $a(b+c)<1$, then the solution $(u(t,x),v(t,x))$ of {\rm(DFB)} satisfies
 \bes
 &\dd\liminf_{t\to\infty}u(t,x)\geq\ud u(x), \ \ \limsup_{t\to\infty}u(t,x)\leq\bar u(x)\ \ \ \mbox{uniformly in any compact subset of } \, [0, \infty),\qquad&\lbl{3.11}\\
 &\dd\liminf_{t\to\infty}v(t,x)\geq\ud v(x), \ \ \limsup_{t\to\infty}v(t,x)\leq\bar v(x)\ \ \ \mbox{uniformly in any compact subset of } \, [0, \infty),\qquad&\lbl{3.12}
 \ees
where $\bar u,\,\bar v,\,\ud u$ and $\ud v$ are given in the proof of Theorem $\ref{th2.3}$.
\end{theo}

{\bf Proof}. \ {\it Step 1} \ Define
 \[\phi(x)=\left\{\begin{array}{ll}
  u_0(x), \ \ & 0\leq x\leq h_0,\\[2mm]
  0, \ \ & x\geq h_0,
  \end{array}\right.\]
and let $w(t,x)$ be the unique positive solution of
 \bess\left\{\begin{array}{lll}
 w_t-w_{xx}=w(1-w), &t>0,\ \ 0<x<\infty,\\[2mm]
 w(t,0)=0,\ \ \ &t>0,\\[2mm]
  w(0,x)=\phi(x), &x\geq 0.
  \end{array}\right.\eess
By the comparison principle, $u(t,x)\leq w(t,x)$ for all $t>0$ and $0\leq x\leq h(t)$.
In view of Theorem \ref{th2.2}, $\lim_{t\to\infty}w(t,x)=\bar u(x)$ uniformly in any compact subset of $[0, \infty)$, where $\bar u(x)$ is the unique positive solution of (\ref{2.16}). Thanks to $h(\infty)=\infty$, we get the second limit of (\ref{3.11}).

{\it Step 2} \ For any given $0<\ep\ll 1$ and $l\gg 1$, there exists a large $T$, such that
 \[h(t)>l, \ \ \ u(t,x)<\bar u(x)+\ep\ff 0\leq x\leq l, \ t\geq T.\]
Let $v^l(t,x)$ be the unique positive solution of
 \bess\left\{\begin{array}{lll}
 v_t-dv_{xx}=v\big[c-v+b(\bar u(x)+\ep)\big],\ \ &t>T, \ \ 0<x<l,\\[2mm]
 v(t,0)=0, \ \ v(t,l)=K,\ \ \ &t>T,\\[2mm]
 v(T,x)=K,&x\in [0,l],
 \end{array}\right.\eess
where $K>\max\big\{M,\, c+b(1+\ep)\big\}$, and $M$ is given by Theorem \ref{th1.1}. Since $v(t,x)\leq M$, by the comparison principle we have
\bes
 v(t,x)\leq v^l(t,x)\ff 0\leq x\leq l, \ t\geq T.\lbl{3.13}\ees
Let $v_\ep(x)$ be the unique positive solution of (\ref{2.2}) with $\lambda=1$ and $f(x)=c+b(\bar u(x)+\ep)$, and let $v_l(x)$  be the unique positive solution of
  \bes\left\{\begin{array}{lll}
 -dv''=v\big[c-v+b(\bar u(x)+\ep)\big],\ \ & 0<x<l,\\[2mm]
 v(0)=0, \ \ v(l)=K.
 \end{array}\right.\lbl{3.14}\ees
Thanks to $v_\ep(x)<K$, the comparison principle asserts
 \[v_\ep(x)\leq v_l(x), \ \ \ v_\ep(x)\leq v^l(t,x) \ff t\geq T, \ \ 0\leq x\leq l.\]
Because $K$ is an upper solution of (\ref{3.14}), it follows that the limit $\lim_{t\to\infty}v^l(t,x)$ exists and is a positive solution of (\ref{3.14}). By the uniqueness of $v_l(x)$ we have that $\lim_{t\to\infty}v^l(t,x)=v_l(x)$ and this limit holds uniformly in $x\in[0,l]$. In view of (\ref{3.13}), it yields
 \bes
  \limsup_{t\to\infty}v(t,x)\leq v_l(x) \ \ \mbox{uniformly on } \ [0,l].\lbl{3.15}\ees
By Proposition \ref{p2.2}, $\lim_{l\to\infty}v_l(x)=v_\ep(x)$ uniformly in any compact
subset of $[0, \infty)$. By Proposition \ref{p2.3}, $\lim_{\ep\to 0}v_\ep(x)=\bar v(x)$ uniformly in any compact subset of $[0, \infty)$, where $\bar v(x)$ is the unique positive solution of (\ref{2.17}). These facts and (\ref{3.15}) imply the second limit of (\ref{3.12}).

{\it Step 3} \ Since $a(b+c)<1$, choose $\ep_0>0$ such that $a(b+c+\ep_0)<1$. For any given $0<\ep<\ep_0$ and $l\gg 1$, there exists a large $T$ such that
 \[h(t)>l, \ \ \ v(t,x)<\bar v(x)+\ep\ff 0\leq x\leq l, \ t\geq T.\]
Moreover, when $l\gg 1$, the problem
  \bess\left\{\begin{array}{lll}
 -u''=u\big[1-u-a(b+c+\ep)\big],\ \ & 0<x<l,\\[2mm]
 u(0)=0=u(l)
  \end{array}\right.\eess
has a unique positive solution, denoted by $u^*_l(x)$. Thanks to $u_x(T, 0)>0$ and $ u(T,l)>0$, there is a positive constant $\sigma<1$ such that $u(T, x)\geq\sigma u^*_l(x)$ for all $0\leq x\leq l$. As $\bar v(x)\leq b+c$, it is easy to see that $\sigma u^*_l(x)$ is a lower solution of the problem
  \bes\left\{\begin{array}{lll}
 -u''=u\big[1-u-a(\bar v(x)+\ep)\big],\ \ & 0<x<l,\\[2mm]
 u(0)=0=u(l).
  \end{array}\right.\lbl{3.16}\ees
Let $u^l(t,x)$ be the unique solution of
 \bess\left\{\begin{array}{lll}
 u_t-u_{xx}=u\big[1-u-a(\bar v(x)+\ep)\big],\ \ &t>T, \ \ 0<x<l,\\[2mm]
 u(t,0)=0=u(t,l),\ \ \ &t>T,\\[2mm]
 u(T,x)=\sigma u_l^*(x),&x\in [0,l].
 \end{array}\right.\eess
Then $u(t,x)\geq u^l(t,x)$ for $0\leq x\leq l$ and $t\geq T$, and $u^l(t,x)$ is increasing in $t$. Similar to the above, it can be deduced that the limit $\lim_{t\to\infty}u^l(t,x):=\ud u_l(x)$ exists and is the unique positive solution of (\ref{3.16}). Moreover, such limit holds uniformly on $[0,l]$. Hence,
 \bes
 \liminf_{t\to\infty}u(t,x)\geq\ud u_l(x) \ \  \mbox{uniformly on} \ \, [0,l].
 \lbl{3.17}\ees
Similarly to the argument of Step 2, let $l\longrightarrow\infty$ firstly and $\ep\longrightarrow 0$ secondly, and apply Propositions \ref{p2.2} and \ref{p2.3} successively. We can prove that
 \[\lim_{l\to\infty}\ud u_l(x)=\ud u(x) \ \ \ \mbox{uniformly in any compact
 subset of } \, [0, \infty).\]
This fact combined with (\ref{3.17}) allows us to derive the first limit of (\ref{3.11}).

Similarly, we can prove the first limit of (\ref{3.12}).
The proof is finished. \ \ \ \fbox{}

\section{Comparison principles}
\setcounter{equation}{0}

In this section we shall provide some comparison principles which will be used to estimate the solution $(u, v, h)$ and determine the criteria governing spreading and vanishing.

\begin{lem} $($Comparison principle$)$\label{l4.1} \
Let $\bar h\in C^1([0,\infty))$ and $\bar h(t)>0$ in $[0,\infty)$. Let $\bar u,\,\bar v\in C(\overline{O})\bigcap C^{1,2}(O)$, with $O=\{(t,x): t>0,\, 0<x<\bar h(t)\}$. Assume that $(\bar u,\bar v, \bar h)$ satisfies
 \bes\left\{\begin{array}{ll}
  \bar u_t-\bar u_{xx}\geq\bar u(1-\bar u),\ \ &t>0,\ \ 0<x<\bar h(t),\\[2mm]
 \bar v_t-d\bar v_{xx}\geq\bar v(c-\bar v+b\bar u),&t>0,\ \ 0<x<\bar h(t),\\[2mm]
 \bar u(t,0)\geq 0,\ \ \bar v(t,0)\geq 0,&t>0,\\[2mm]
 \bar u(t,\bar h(t))=\bar v(t,\bar h(t))=0,\ \ &t>0,\\[2mm]
 \bar h'(t)\geq-\mu\big[\bar u_x(t,\bar h(t))
 +\rho\bar v_x(t,\bar h(t))\big],\ \ &t>0.
 \end{array}\right.\lbl{4.1}\ees
If $\bar h(0)\geq h_0$, $\bar u(0,x),\, \bar v(0,x)\geq 0$ on $[0,\bar h(0)]$, and $\bar u(0,x)\geq u_0(x),\,\bar v(0,x)\geq v_0(x)$ on $[0,h_0]$. Then the solution $(u,v,h)$ of {\rm(DFB)} satisfies $h(t)\leq\bar h(t)$ on $[0,\infty)$, and $u\leq\bar u,\,v\leq\bar v$ on $D$, where $D=\{(t,x): t\geq 0,\, 0\leq x\leq h(t)\}$.

If, in $(\ref{4.1})$, the conditions $\bar u(t,0)\geq 0$ and $\bar v(t,0)\geq 0$ are replaced by
$\bar u_x(t,0)\leq 0$ and $\bar v_x(t,0)\leq 0$, then the conclusion still holds for the solution of {\rm(NFB)}.
\end{lem}

{\bf Proof}. \ The proof can proceed as the argument of \cite[Lemma 5.1]{GW} with minor modification.  We first consider that $\bar h(0)>h_0$. Then $\bar h(t)>h(t)$ for small $t>0$. We can derive that $\bar h(t)>h(t)$ for all $t\geq 0$. If this is not true, there exists $t_0>0$ such that $\bar h(t_0)=h(t_0)$ and $\bar h(t)>h(t)$ for all $t\in(0,t_0)$. Thus, $\bar h'(t_0)\leq h'(t_0)$. Set
 \[D_{t_0}=\bbb\{(t,x):\,t\in(0,t_0],\ 0<x<h(t)\bbb\}.\]
Recall that $\bar u(0,x)\geq u_0(x)$ on $[0,h_0]$ and $\bar u(t_0, h(t_0))=0=u(t_0,h(t_0))$, the strong maximal principle yields that $\bar u>u$ in $D_{t_0}$, and $\bar u_x(t_0,h(t_0))<u_x(t_0, h(t_0))$. Similarly, $\bar v>v$ in $D_{t_0}$, and $\bar v_x(t_0,h(t_0))<v_x(t_0, h(t_0))$. However,
 \[\bar h'(t_0)\geq-\mu\big[\bar u_x(t_0,h(t_0))+\rho\bar v_x(t_0,h(t_0))>
 -\mu\big[u_x(t_0,h(t_0))+\rho v_x(t_0,h(t_0))\big]=h'(t_0).\]
We get a contradiction. Hence, $\bar h(t)>h(t)$ for all $t\geq 0$, and $u\leq\bar u,\,v\leq\bar v$ on $\overline{D}$.

When $\bar h(0)=h_0$, the process is the same as that of \cite[Lemma 5.1]{GW}. \ \ \ \fbox{}

In the same way as the proof of Lemma \ref{l4.1} we can prove the following lemma.

\begin{lem} $($Comparison principle$)$\label{l4.2}\
Let $\ud h\in C^1([0,\infty))$ with $\ud h(t)>0$ for
all $t\in [0,\infty)$, and $\ud u\in C(\overline{O}_1)\bigcap C^{1,2}(O_1)$ with
$O_1=\{(t,x): t>0,\, 0<x<\ud h(t)\}$.
Assume that $(\ud v, \ud h)$ satisfies
 \bes
 \left\{\begin{array}{ll}
  \ud v_t-d\ud v_{xx}\leq\ud v(c-\ud v),\ \ &t>0,\, \ 0<x<\ud h(t),\\[2mm]
  \ud v(t,0)=\ud v(t,\ud h(t))=0,\ \ &t>0,\\[2mm]
 \ud h'(t)\leq-\mu\rho\ud v_x(t,\ud h(t)),\ \ &t>0
 \end{array}\right.\lbl{4.2}
 \ees
and $0<\ud h(0)\leq h_0$, $0\leq\ud v(0,x)\leq v_0(x)$ on $[0,\ud h(0)]$.
Then the solution $(u,v,h)$ of {\rm(DFB)} satisfies $h(t)\geq\ud h(t)$ on $[0,\infty)$, and $v(t,x)\geq\ud v(t,x)$ on $\overline{O}_1$.

If, in $(\ref{4.2})$, the condition $\ud v(t,0)=0$ is replaced by
$\ud v_x(t,0)\geq 0$, then the conclusion still holds for the solution of {\rm(NFB)}.
\end{lem}



\section{The criteria governing spreading and vanishing}
\setcounter{equation}{0}

In this section we study the criteria governing spreading and vanishing for the problems (DFB) and (NFB), respectively. We first give a necessary condition of vanishing.

\begin{theo}\lbl{th5.1} \ If $h_\infty<\infty$, then $h_\infty\leq\pi\min\big\{\sqrt{d/c},\,1\big\}$ for the problem {\rm(DFB)}, and $h_\infty\leq\frac {\pi}2\min\big\{\sqrt{d/c},\,1\big\}$ for the problem {\rm(NFB)}.

Define
  \[\Lambda=\lk\{\begin{array}{ll}\medskip
  \pi\min\big\{\sqrt{d/c},\,1\big\} \ \ \mbox{ for the problem} \ {\rm(DFB)},\\
  \dd\frac{\pi}2\min\big\{\sqrt{d/c},\,1\big\} \ \ \mbox{ for the problem} \ {\rm(NFB)}.
  \end{array}\right.\]
Then $h_0\geq \Lambda$ implies $h_\infty=\infty$ due to $h'(t)>0$ for $t>0$.
\end{theo}

{\bf Proof}. We only deal with the problem (DFB), since the problem (NFB) can be treated by the similar way. By Theorem \ref{th3.2}, $h_\infty<\infty$ implies
$\lim\limits_{t\to\infty}\|u(t,\cdot), v(t,\cdot)\|_{C[0,\,h(t)]}=0$. We assume $h_\infty>\Lambda$ to get a contradiction.

If $h_\infty>\pi$, then there exists $\ep>0$ such that $h_\infty>\pi\sqrt{1/(1-a\ep)}$. For such $\varepsilon$, there exists $T\gg 1$ such that $h(T)>\pi\sqrt{1/(1-a\ep)}$ and
  \[v(t,x)\leq \varepsilon,\ \ \ \forall \ t\geq T, \ x\in[0,h(T)].\]
Set $l=h(T)$ and let $w=w(t,x)$ be the unique positive solution of the following initial boundary value problem with fixed boundary:
  $$\left\{\begin{array}{ll}
 w_t=w_{xx}+w\dd\lk(1-w-a\varepsilon\right), \ \ &t>T, \ \, 0<x<l,\\[2mm]
 w(t,0)=w(t,l)=0,&t>T,\\[2mm]
  w(T,x)=u(T,x),&0\leq x\leq l.
  \end{array}\right.$$
By the comparison principle,
 $$w(t,x)\leq u(t,x)\ff t\geq T,\ 0\leq x\leq l.$$
Since $l>\pi\sqrt{1/(1-a\ep)}$,
it is well known that $w(t,x)\longrightarrow W(x)$ as $t\longrightarrow\infty$ uniformly in any  compact subset of $(0,l)$, where $W$ is the unique positive solution of
  \[\left\{\begin{array}{ll}
 W_{xx}+W\dd\lk(1-a\varepsilon-W\right)=0,\ \ &0<x<l,\\[2mm]
 W(0)=W(l)=0.&\end{array}\right.\]
Hence, $\dd\liminf_{t\to\infty} u(t,x)\geq\lim_{t\to\infty}w(t,x)=W(x)>0$ in
$(0,l)$. This is a contradiction to (\ref{3.10}).

\vskip 4pt If $h_\infty>\pi\sqrt{d/c}$, then there exists $T\gg 1$ such that $h(T)>\pi\sqrt{d/c}$. Set $l=h(T)$ and let $z=z(t,x)$ be the unique positive solution of the following initial boundary value problem with fixed boundary:
  $$\left\{\begin{array}{ll}
 z_t=dz_{xx}+z\lk(c-z\right), \ \ &t>T, \ \, 0<x<l,\\[2mm]
 z(t,0)=z(t,l)=0,&t>T,\\[2mm]
  z(T,x)=v(T,x),&0\leq x\leq l.
  \end{array}\right.$$
By the comparison principle,
 $$z(t,x)\leq v(t,x)\ff t\geq T,\ 0\leq x\leq l.$$
Since $l>\pi\sqrt{d/c}$, similarly to the above, we can get a contradiction to (\ref{3.10}).
\ \ $\Box$

\vskip 2pt Now we discuss the case $h_0<\Lambda$.

\begin{lem}\lbl{l5.1} \ Suppose that $h_0<\Lambda$. For the problem {\rm(DFB)}, if
  \bess
  \mu\geq\mu^0:=\max\lk\{1,\dd\frac{1}{c}\|v_0\|_\infty\right\}\frac d\rho\lk(\frac{\pi^2d}{c}
  -h^2_0\right)\lk(2\dd\int_0^{h_0}xv_0(x)dx\right)^{-1},
  \eess
then $h_\infty=\infty$. For the problem {\rm(NFB)}, if
  \bess
  \mu\geq\mu^0:=\max\lk\{1,\dd\frac{1}{c}\|v_0\|_\infty\right\}\frac d\rho\lk(\frac{\pi}2\sqrt{\frac cd}-h_0\right)\lk(\dd\int_0^{h_0}v_0(x)dx\right)^{-1},
  \eess
then $h_\infty=\infty$.
\end{lem}

{\bf Proof}. \ For the problem (DFB), we consider the following auxiliary problem
$$
 \left\{\begin{array}{ll}
  \underline v_t-d\underline v_{xx}=\underline v(c-\underline v),
  \ &t>0, \ \ 0<x<\underline h(t),\\[2mm]
 \underline v(t,0)=0,\ \ \underline v(t,\underline h(t))=0,\ &t>0,\\[2mm]
 \underline h'(t)=-\mu\rho\underline v_x(t,\underline h(t)), \ &t>0,\\[2mm]
 \underline v(0,x)=v_0(x),\ \ \ &0\leq x\leq h_0,\\[2mm]
 \underline h(0)=h_0.
 \end{array}\right.$$
It follows from Lemma \ref{l4.2} that
  $$\underline h(t)\leq h(t),\ \ \underline v(t,x)\leq v(t,x)
  \ff t>0,\ \ 0<x<\underline h(t).$$
Recall that $h_0<\Lambda\leq\pi\sqrt{d/c}$, and $\mu\geq\mu^0$, in view of the Proposition 4.8 in \cite{KY}, it yields $\underline
h(\infty)=\infty$. Therefore, $h_\infty=\infty$.

For the problem (NFB), we consider the following auxiliary problem
$$
 \left\{\begin{array}{ll}
  \underline v_t-d\underline v_{xx}=\underline v(c-\underline v),
  \ &t>0, \ \ 0<x<\underline h(t),\\[2mm]
 \underline v_x(t,0)=0,\ \ \underline v(t,\underline h(t))=0,\ &t>0,\\[2mm]
 \underline h'(t)=-\mu\rho\underline v_x(t,\underline h(t)), \ &t>0,\\[2mm]
 \underline v(0,x)=v_0(x),\ \ \ &0\leq x\leq h_0,\\[2mm]
 \underline h(0)=h_0.
 \end{array}\right.$$
By Lemma \ref{l4.2},
  $$\underline h(t)\leq h(t),\ \ \underline v(t,x)\leq v(t,x)
  \ff t>0, \,\ 0<x<\underline h(t).$$
Note that $h_0<\Lambda\leq\frac{\pi}2\sqrt{d/c}$, and $\mu\geq\mu^0$, by use of the Lemma 3.7 of \cite{DLin}, we have $\underline
h(\infty)=\infty$. Therefore, $h_\infty=\infty$. \ \ \ \fbox{}

\begin{lem}\lbl{l5.2} \ Assume that $h_0<\Lambda$. Then
there exists $\mu_0>0$, depending also on $u_0(x)$ and $v_0(x)$, such that
$h_\infty<\infty$ when $\mu\leq\mu_0$ for both problems {\rm(DFB)} and {\rm(NFB)}.
\end{lem}

{\bf Proof}. \ We shall use the argument from Ricci and Tarzia \cite{RT} to construct suitable upper solutions and use Lemma \ref{l4.1} to derive the desired conclusion.

{\it Step 1} \, We first consider the problem (NFB) and adopt the following
functions constructed by Du and Lin \cite{DLin} (see also Guo and Wu \cite{GW}):
 \bess
 &\dd\sigma(t)=h_0\lk(1+\delta-\frac\delta 2{\rm e}^{-\beta t}\right), \ \ t\geq 0; \quad \
 V(y)=\cos\dd\frac{\pi y}{2},\ \ 0\leq y\leq 1,& \\[1mm]
 &\bar u(t,x)=Me^{-\beta t}V\lk(\dd\frac x{\sigma(t)}\right), \ \
 \bar v(t,x)=bMe^{-\beta t}V\lk(\dd\frac x{\sigma(t)}\right), \ \ 0\leq x\leq\sigma(t).&
 \eess
It is obvious that
 \[\bar u_x(t,0)=\bar v_x(t,0)=0, \ \ \bar u\bbb(t,\sigma(t)\bbb)
 =\bar v\bbb(t,\sigma(t)\bbb)=0\ff t\geq 0.\]
Recall that $h_0<\frac{\pi}2\min\big\{\sqrt{d/c},\,1\big\}$ in the present case and $\bar v(t,x)=b\bar u(t,x)$. Similarly to the proof of Corollary 1(iii) in \cite{GW} (pp.892), we can verify that, for suitable small positive constants $\delta$ and $\beta$, and large positive constant $M$, the pair $(\bar u, \bar v)$ satisfies
 \bes\left\{\begin{array}{ll}\medskip
\bar u_t-\bar u_{xx}-\bar u(1-\bar u)\geq 0,\ \ &t>0, \ \ 0<x\leq\sigma(t),\\
\bar v_t-d\bar v_{xx}-\bar v(c-\bar v+b\bar u)\geq 0,\ \ &t>0, \ \ 0<x\leq\sigma(t),\\[2mm]
\bar u(0,x)\geq u_0(x),\ \ \bar v(0,x)\geq v_0(x),\ \ &0\leq x\leq h_0(1+\delta/2).
\nm\end{array}\right.\ees
Moreover, for such fixed constants $\delta$, $\beta$ and $M$, there exists $\mu_0>0$ such that
 \[\sigma'(t)+\mu\big(\bar u_x+\rho\bar v_x\big)\big|_{x=\sigma(t)}\geq 0\]
for all $\mu\leq\mu_0$.

By Lemma \ref{l4.1}, $\sigma(t)\geq h(t)$. Taking $t\longrightarrow\infty$ we have $h_\infty\leq\sigma(\infty)=h_0(1+\delta)<\infty$.

\vskip 4pt {\it Step 2} \, Now we discuss the problem (DFB). Recall that $h_0<\pi\min\big\{\sqrt{d/c},\,1\big\}$ for our present case. We can verify that there exist two positive constants $\delta,\beta\ll 1$ such that
 \bes
 \frac 1{h_0(1+\delta)}\lk(\frac{\pi^2}{h_0(1+\delta)}-\beta\delta h_0\right)-\beta-1>0,\lbl{q.1}\\
 \frac 1{h_0(1+\delta)}\lk(\frac{d\pi^2}{h_0(1+\delta)}-\beta\delta h_0\right)-\beta-c>0.\lbl{q.2}\ees
For such fixed $\delta$ and $\beta$, let
 \bess
 &\dd\sigma(t)=h_0\lk(1+\delta-\frac\delta 2{\rm e}^{-\beta t}\right), \ \ t\geq 0; \quad \
 W(y)=\sin\pi y,\ \ 0\leq y\leq 1,& \\
 &\bar u(t,x)=Me^{-\beta t}W\lk(\dd\frac x{\sigma(t)}\right), \ \
 \bar v(t,x)=bMe^{-\beta t}W\lk(\dd\frac x{\sigma(t)}\right), \ \ 0\leq x\leq\sigma(t),&
 \eess
where $M$ is a large positive constant. It is obvious that
 \bes\bar u(t,0)=\bar v(t,0)=\bar u\bbb(t,\sigma(t)\bbb)
 =\bar v\bbb(t,\sigma(t)\bbb)=0\ff t\geq 0,\lbl{5.1}\ees
and
 \bes\bar u(0,x)\geq u_0(x),\ \ \bar v(0,x)\geq v_0(x)\ff 0\leq x\leq\sigma_0(1+\delta/2)\lbl{5.2}\ees
provided that $M\gg 1$. Moreover, for such fixed constants $\delta$, $\beta$ and $M$, since $\sigma(t)\geq h_0(1+\delta/2)$, it is easy to see that there exists $0<\mu_0\ll 1$ such that
 \bes
 \sigma'(t)+\mu\big(\bar u_x+\rho\bar v_x\big)\big|_{x=\sigma(t)}={\rm e}^{-\beta t}\lk(\frac{\beta\delta h_0}2-\mu M\pi(1+b\rho)\frac 1{\sigma(t)}\right)>0\ff t\geq 0\lbl{5.3}\ees
provided that $0<\mu\leq\mu_0$.

Denote $y=x/\sigma(t)$. The direct calculation yields,
 \bess
 \bar u_t-\bar u_{xx}-\bar u(1-\bar u)&=&\bar u\lk(\frac{\pi^2}{\sigma^2(t)}-\beta-1
 -{\rm e}^{-\beta t}\frac{\beta\delta h_0}{\sigma(t)}\times\frac{\pi y}{2\sin\pi y}\cos\pi y+\bar u\right)\\
 &\geq&\bar u\lk(\frac{\pi^2}{\sigma^2(t)}-\beta-1
 -{\rm e}^{-\beta t}\frac{\beta\delta h_0}{\sigma(t)}\times\frac{\pi y}{2\sin\pi y}\cos\pi y\right).
 \eess
Since $\cos\pi y\leq 0$ for $1/2\leq y\leq 1$, and $\sigma(t)$ is increasing, we have
 \bess
 \bar u_t-\bar u_{xx}-\bar u(1-\bar u)
 &\geq&\bar u\lk(\frac{\pi^2}{\sigma^2(t)}-\beta-1\right)\\[1mm]
 &\geq&\bar u\lk(\frac{\pi^2}{h_0^2(1+\delta)^2}-\beta-1\right)\\[2mm]
 &>&0\ff \sigma(t)/2\leq x\leq\sigma(t)
 \eess
by (\ref{q.1}). Remember that $0\leq\cos\pi y\leq 1$, $y\leq\frac2\pi\sin\pi y$ for all $0\leq y\leq 1/2$, and ${\rm e}^{-\beta t}\leq 1$ for all $t\geq 0$. We have that, for all $t>0$ and $0\leq x\leq\sigma(t)/2$,
 \[{\rm e}^{-\beta t}\frac{\beta\delta h_0}{\sigma(t)}\times\frac{\pi y}{2\sin\pi y}\cos\pi y
 \leq \frac{\beta\delta h_0}{\sigma(t)}.\]
It follows that, for all $t>0$ and $0\leq x\leq\frac 12\sigma(t)$,
 \bess
 \bar u_t-\bar u_{xx}-\bar u(1-\bar u)&\geq&\bar u\lk(\frac{\pi^2}{\sigma^2(t)}-\beta-1
 -\frac{\beta\delta h_0}{\sigma(t)}\right)\\
 &=&\bar u\lk[\frac1{\sigma(t)}\lk(\frac{\pi^2}{\sigma(t)}-\beta\delta h_0\right)-\beta-1\right]\\
 &\geq &\bar u\lk[\frac 1{h_0(1+\delta)}\lk(\frac{\pi^2}{h_0(1+\delta)}-\beta\delta h_0\right)-\beta-1\right]\\[1mm]
 &>&0
 \eess
by (\ref{q.1}). In conclusion, we have
 \bes\bar u_t-\bar u_{xx}-\bar u(1-\bar u)>0 \ff t>0, \ \ 0\leq x\leq\sigma(t).\lbl{5.4}\ees

Remember $\bar v(t,x)=b\bar u(t,x)$; in view of (\ref{q.2}), similar to the above, we can verify that
  \bes
  \bar v_t-d\bar v_{xx}-\bar v(c-\bar v+b\bar u)>0 \ff t>0, \ \ 0\leq x\leq\sigma(t).\lbl{5.5}\ees

Notice that (\ref{5.1})--(\ref{5.5}), by virtue of Lemma \ref{l4.1} we have $\sigma(t)\geq h(t)$. Taking $t\longrightarrow\infty$ we have $h_\infty\leq\sigma(\infty)=h_0(1+\delta)<\infty$.
The proof is complete.  \ \ \ \fbox{}

\begin{theo}\lbl{th5.2}\, Suppose that $h_0<\Lambda$. Then there exist $\mu^*\geq\mu_*>0$, depending on $u_0(x)$, $v_0(x)$ and $h_0$, such that $h_\infty=\infty$ if $\mu>\mu^*$, and $h_\infty\leq\Lambda$ if $\mu\leq\mu_*$ or $\mu=\mu^*$.
\end{theo}

{\bf Proof}.\, The proof is similar to that of \cite[Theorem 3.9]{DLin} and \cite[Theorem 5.4]{WZ}. For the convenience to the reader we shall give the details
because this is a main theorem in this section. We will
write $(u_\mu, v_\mu, h_\mu)$ in place of $(u, v, h)$ to clarify the
dependence of the solution of (DFB) and/or (NFB) on $\mu$. Define
 $$\Sigma^*=\lk\{\mu>0:\,h_{\mu,\infty}\leq\Lambda\right\}.$$
By Lemma \ref{l5.2}, $(0,\mu_0]\subset\Sigma^*$. In view of Lemma \ref{l5.1}, $\Sigma^*\bigcap[\mu^0,\infty)=\emptyset$. Therefore,
$\mu^*:=\sup\Sigma^*\in[\mu_0,\,\mu^0]$. By this definition and Theorem \ref{th5.1} we
find that $h_{\mu,\infty}=\infty$ when $\mu>\mu^*$.
Hence, $\Sigma^*\subset(0,\mu^*]$.

We will show that $\mu^*\in\Sigma^*$. Otherwise, $h_{\mu^*,\infty}=\infty$. There exists $T>0$ such that $h_{\mu^*}(T)>\Lambda$. By the continuous
dependence of $(u_\mu, v_\mu, h_\mu)$ on $\mu$, there is $\varepsilon>0$ such
that $h_{\mu}(T)>\Lambda$ for
$\mu\in(\mu^*-\varepsilon,\mu^*+\varepsilon)$. It follows that for all such $\mu$,
 \[\lim_{t\to\infty}h_{\mu}(t)\geq h_{\mu}(T)>\Lambda.\]
Therefore, $\bbb(\mu^*-\varepsilon,\mu^*+\epsilon\bbb)\bigcap\Sigma^*=\emptyset$, and
$\sup\Sigma^*\leq\mu^*-\varepsilon$. This contradicts the definition of $\mu^*$.

Define
 $$\Sigma_*=\lk\{\nu:\, \nu\geq\mu_0 \ \mbox{such\ that} \  \,h_{\mu,\infty}\leq\Lambda\ \mbox{for\ all} \ 0<\mu\leq\nu\right\},$$
where $\mu_0$ is given by Lemma \ref{l5.2}. Then $\mu_*:=\sup\Sigma_*\leq\mu^*$ and $(0,\mu_*)\subset\Sigma_*$. Similarly to the above, we can prove that $\mu_*\in\Sigma_*$. The proof is completed.
 \ \ \ \fbox{}

\section{Asymptotic spreading speed}
\setcounter{equation}{0}

In this section we provide an upper bound for $\dd\limsup_{t\to\infty}\frac{h(t)}t$, which shows that the asymptotic spreading speed (if exists) for both problems (DFB) and (NFB) cannot be faster than $2\max\big\{\sqrt{cd},\,1\big\}$ under some suitable conditions. The number $2\max\big\{\sqrt{cd},\,1\big\}$ seems to be the minimal speed of traveling wave fronts of the prey-predator system
 \bes\left\{\begin{array}{lll}
 u_t-u_{xx}=u(1-u-av), &t>0,\ \ x\in\mathbb{R},\\[2mm]
  v_t-dv_{xx}=v(c-v+bu),\ \ &t>0, \ \ x\in\mathbb{R},
  \end{array}\right.\lbl{6.1}\ees
please refer to \cite{LinG}.

It is easy to see that Theorem 5.17 of \cite{LLM} still holds for the traveling wave fronts of the following system
 $$
 \left\{\begin{array}{lll}
 u_t-u_{xx}=u(1-u), &t>0,\ \ x\in\mathbb{R},\\[2mm]
  v_t-dv_{xx}=v(c-v+bu),\ \ &t>0, \ \ x\in\mathbb{R}.
 \end{array}\right.$$
Thus, for any given $s>2\max\big\{\sqrt{cd},\,1\big\}$, the following problem
 \bes\left\{\begin{array}{ll}\medskip
 s\phi'+ \phi''+ \phi(1-\phi)=0, \quad
 s\psi'+d\psi''+\psi(c-\psi+b\phi)=0 \ \ {\rm in} \ \ \mathbb{R},\\
 \phi(-\infty)=1, \ \ \psi(-\infty)=b+c, \ \ (\phi,\psi)(\infty)=(0,0), \ \ \phi(0)=1/2,\\[2mm]
\phi'<0, \ \ \psi'<0  \ \ {\rm in} \ \ \mathbb{R}\lbl{6.2}\end{array}\right.\ees
has a solution $(\phi(\xi),\psi(\xi))$, with $\xi=x-st$. Moreover, $(\phi(\xi),\psi(\xi))$ satisfies
  \bes\lim_{\xi\to\infty}\phi(\xi){\rm e}^{\lambda_1\xi}=\lim_{\xi\to\infty}\psi(\xi){\rm e}^{\lambda_2\xi}=1,
  \lbl{6.3}\ees
where
 \bes
 \lambda_1=\frac{s+\sqrt{s^2-4}}2>0, \ \ \ \lambda_2=\frac{s+\sqrt{s^2-4cd}}{2d}>0.
 \lbl{6.4}\ees

\begin{theo}\lbl{th6.1}\, Let $(u,v,h)$ be the solution of the problem {\rm(DFB)} or {\rm(NFB)} and $h_\infty=\infty$. If for any given $s>2\max\big\{\sqrt{cd},\,1\big\}$, the problem {\rm(\ref{6.2})} has a solution $(\phi(\xi),\psi(\xi))$ satisfying
 \bes
 \psi(\xi)\geq\beta\phi(\xi)\ff \xi\in\mathbb{R}\lbl{6.5}\ees
for some positive constant $\beta$ $($which may depend on $s)$. Then we have
 \bess
 \limsup_{t\to\infty}\frac{h(t)}t\leq 2\max\big\{\sqrt{cd},\,1\big\}.
 \eess
 \end{theo}

Before giving the proof of Theorem \ref{th6.1}, we state one remark to guarantee the condition (\ref{6.5}).

\begin{remark}\lbl{r6.1}\, For any given $s>2\max\big\{\sqrt{cd},\,1\big\}$, let $\lambda_1$ and $\lambda_2$ be given by {\rm(\ref{6.4})}. By the carefully calculations we have that if one of the following holds:
 \vspace{-2mm}\begin{quote}
{\rm(a)}\, $d\geq 1$, $cd\geq 1$; \\
{\rm(b)}\, $d\geq 1$, $cd<1$ and $c+d\geq 2$; \\
{\rm(c)}\, $d<1$, $cd\geq 1$ and $2cd\geq c+d$,
 \vspace{-2mm}\end{quote}
then $\lambda_2\geq\lambda_1$. In view of {\rm(\ref{6.3})} and the limits:
 \[\lim_{\xi\to-\infty}\phi(\xi)=1, \ \ \lim_{\xi\to-\infty}\psi(\xi)=b+c,\]
it can be seen that there exists a positive constant $\beta$ such that {\rm(\ref{6.5})} holds. \end{remark}

{\bf Proof of Theorem \ref{th6.1}}.\, The idea of this proof comes from \cite{GW}. For any given $s>2\max\big\{\sqrt{cd},\,1\big\}$, let $(\phi(\xi),\,\psi(\xi))$ be the solution of {\rm(\ref{6.2})} satisfying {\rm(\ref{6.5})}.

Choose $m>k\gg 1$ such that
 \bes
 m\beta \geq 2bk, \ \ \ \mbox{which implies} \ \ (m-1)\beta \geq 2b(k-1), \lbl{6.6}\\
 k\phi(\xi)>\|u_0\|_\infty, \ \ \ m\psi(\xi)>\|v_0\|_\infty\ff\xi\in[0,h_0].\nm\ees
For such fixed $m$ and $k$, recall that $(\phi(\xi),\psi(\xi))\longrightarrow 0$
and $(\phi'(\xi),\psi'(\xi))\longrightarrow 0$ as $\xi\longrightarrow\infty$, there exists $\sigma_0>h_0$ such that
 \bes
 &\dd \phi(\sigma_0)<\min_{0\leq x\leq h_0}\lk(\phi(x)-\frac {u_0(x)}k\right), \ \
 \psi(\sigma_0)<\min_{0\leq x\leq h_0}\lk(\psi(x)-\frac {v_0(x)}m\right),&\lbl{6.7}\\[2mm]
 &\phi(\sigma_0)<1-1/k, \ \ \psi(\sigma_0)<2c(m-1)/(m+3m^2),&\lbl{6.8}\\[2mm]
 &-\mu\big[k\phi'(\sigma_0)+m\rho \psi'(\sigma_0)\big]<s.&\lbl{6.9}\ees

Set $\sigma(t)=\sigma_0+st$ and
 \[\bar u(t,x)=k\phi(x-st)-k\phi(\sigma_0), \ \  \bar v(t,x)=m\psi(x-st)-m\psi(\sigma_0).\]
It is obvious that
 \[\bar u(t,\sigma(t))=\bar v(t,\sigma(t))=0 \ff t\geq 0.\]
Since $\phi'<0$, $\psi'<0$, we see that
 \[\bar u(t,0)>0, \ \ \bar v(t,0)>0, \ \ \bar u_x(t,0)<0, \ \ \bar v_x(t,0)<0\ff t\geq 0.\]
It is deduced from (\ref{6.7}) that
 \[\bar u(0,x)>u_0(x), \ \ \bar v(0,x)>v_0(x)\ff 0\leq x\leq h_0.\]
By the first inequality of (\ref{6.8}),
 \bess
 \bar u_t-\bar u_{xx}-\bar u(1-\bar u)=k\lk[(k-1)\lk(\phi-\frac{k \phi(\sigma_0)}{k-1}\right)^2+ \phi(\sigma_0)\frac{k-1-k\phi(\sigma_0)}{k-1}\right]\geq 0.\eess
In order to save space, we denote $\phi(\sigma_0)$ and $\psi(\sigma_0)$ by $\phi_0$ and $\psi_0$, respectively. Applying the second inequality of (\ref{6.8}), (\ref{6.5}) and (\ref{6.6}) we have
 \bess
 &&\frac 1m\big[\bar v_t-d\bar v_{xx}-\bar v(c-\bar v+b\bar u)\big]\\[1mm]
 &=&\frac{(m-1)}2\lk(\psi-\frac{2m}{m-1} \psi_0\right)^2+\psi_0\lk(c-\frac{m(1+3m)}{2(m-1)}\psi_0\right)
 +\frac 12\psi\big[(m-1)\psi-2b(k-1)\phi\big]\\[1mm]
 &&+\frac 12\psi_0\big(m\psi_0-2bk\phi_0\big)+bk(\phi\psi_0+\phi_0\psi)\\[1mm]
 &>&\frac 12\psi\phi\big[(m-1)\beta -2b(k-1)\big]+\frac 12\psi_0\phi_0\big(m\beta -2bk\big)\\[1mm]
 &\geq&0.
   \eess
It follows from (\ref{6.9}) that
 \[\sigma'(t)=s>-\mu \lk[k\phi'(\sigma_0)+\rho m\psi'(\sigma_0)\right]=-\mu\big[\bar u_x(t,\sigma(t))+\rho\bar v_x(t,\sigma(t))\big].\]
Recall $\sigma_0>h_0$; so we have verified the conditions of Lemma \ref{l4.1}.

In view of Lemma \ref{l4.1}, $\sigma(t)\geq h(t)$. Therefore
 \[\limsup_{t\to\infty}\frac{h(t)}t\leq \lim_{t\to\infty}\frac{\sigma(t)}t=s.\]
By the arbitrariness of $s>2\max\big\{\sqrt{cd},\,1\big\}$ we get
 \[\limsup_{t\to\infty}\frac{h(t)}t\leq 2\max\big\{\sqrt{cd},\,1\big\}.\]
The proof is complete.  \ \ \ \fbox{}

\begin{remark}\lbl{r6.4} \, Let $(w,g)$ and $(z,p)$ be solutions of the free boundary problems
 $$
 \left\{\begin{array}{lll}
 w_t-w_{xx}=w(1-w), &t>0,\ \ 0<x<g(t),\\[2mm]
 w_x(t,0)=w(t,g(t)=0,\ \ \ &t>0,\\[2mm]
g'(t)=-\eta w_x(t,g(t)), \ &t>0,\\[2mm]
  w(0,x)=w_0(x), &x\in [0,g_0],\\[2mm]
 g(0)=g_0\geq h_0
 \end{array}\right.$$
and
 $$
 \left\{\begin{array}{lll}
 z_t-dz_{xx}=z(c-z),\ \ &t>0, \ \ 0<x<p(t),\\[2mm]
 z_x(t,0)=z(t,p(t))=0,\ \ \ &t>0,\\[2mm]
 \eta'(t)=-\zeta z_x(t,p(t)), \ &t>0,\\[2mm]
  z(0,x)=z_0(x), &x\in [0,p_0],\\[2mm]
 p(0)=p_0\leq h_0,
 \end{array}\right. $$
respectively. Assume that $g(\infty)=\infty$ and $p(\infty)=\infty$ $($these will be true  under the suitable conditions on the parameters $g_0$, $p_0$, $\eta$ and $\zeta$, refer to {\rm\cite{DLin})}. It follows from the result of {\rm\cite{DLin}} that there are positive constants $g^*$ and $p^*$, such that
 \[\lim_{t\to\infty}\frac{g(t)}t=g^*, \quad \lim_{t\to\infty}\frac{p(t)}t=p^*.\]
Suppose that $\eta\geq\mu$, $\zeta\leq\mu\rho$, and
 \[u_0(x)\leq w_0(x) \ \ \mbox{in} \ \ [0,h_0], \qquad
  v_0(x)\geq z_0(x)>0 \ \ \mbox{in} \ \ [0,p_0].\]
For the solution $(u,v,h)$ of the problem {\rm(NFB)}, by the comparison principle $($Lemma $\ref{l4.2}$ and its similar version$)$ we have that $p(t)\leq h(t)\leq g(t)$. Hence,
 \[p^*\leq\liminf_{t\to\infty}\frac{h(t)}t, \ \ \limsup_{t\to\infty}\frac{h(t)}t\leq g^*.\]
\end{remark}

\section{Discussion}
\setcounter{equation}{0}

In this paper, we have examined a Lotka-Volterra type prey-predator model with free boundary
$x=h(t)$ for both prey and predator, which describes the movement process through
the free boundary.  We envision that the two species initially occupy the region
$[0, h_0]$ and have a tendency to expand their territory together. Then we extend some results of \cite{DLin} and \cite{KY} for one species case and \cite{GW} for two-species weak competition system case to the prey-predator system. The dynamic behavior are discussed. Let $\Lambda=\pi\min\big\{\sqrt{d/c},\,1\big\}$ for the problem (DFB), and
$\Lambda=\frac{\pi}2\min\big\{\sqrt{d/c},\,1\big\}$ for the problem (NFB).
It was proved that:

(i) \, If the size of initial habitat is not less than $\Lambda$, or it is less
than $\Lambda$ but the moving parameter/coefficient $\mu$ of the free boundary is
greater than $\mu^*$ (it depends on the initial data $(u_0,v_0)$ and $h_0$), then
the two species will spread successfully. Moreover,

(ia)\, To the problem (NFB), as $t\longrightarrow\infty$, both $u(t,x)$ and $v(t,x)$ go to positive
constants for the weakly hunting case: $ac<1$ and $ab<1$; while
$u(t,x)\longrightarrow 0$ and $v(t,x)\longrightarrow c$ for the
strongly hunting case: $ac>1$;

\vskip 2pt (ib)\, To the problem (DFB), if $a(b+c)<1$, then $u(t,x)$ and $v(t,x)$ satisfy
\[\liminf_{t\to\infty}u(t,x)\geq\ud u(x), \ \ \limsup_{t\to\infty}u(t,x)\leq\bar u(x),
\ \ \ \liminf_{t\to\infty}v(t,x)\geq\ud v(x), \ \ \limsup_{t\to\infty}v(t,x)\leq\bar v(x)\]
uniformly in any compact subset of $[0, \infty)$, where $\bar u,\,\bar v,\,\ud u$ and
$\ud v$ are positive functions given by Theorem $\ref{th2.3}$.

(ii) \, While if the size of initial habitat is less than $\Lambda$ and the moving parameter/coefficient $\mu$ of the free boundary is less than $\mu_*$, then  $\lim_{t\to\infty}h(t)<\infty$, and $\|u(t,x), v(t,x)\|_{C[g(t),h(t)]}\to 0$ as $t\longrightarrow\infty$. That is, the two species will disappear eventually.

The above conclusions not only provide the spreading-vanishing dichotomy and criteria governing spreading and vanishing, but also provided the long time behavior of $(u(t,x),\,v(t,x))$. If the size of initial habitat is small, and the moving parameter is small enough, it causes no population can survive eventually, while they can coexist if the size of habitat or the moving parameter is large enough, regardless of initial population size. This phenomenon suggests that the size of the initial habitat and the moving parameter are important to the survival for the two species. It is well-known that the effect of habitat size to the survival for species with Dirichlet boundary problem is quite important (see, for example, \cite{BDK}).

Finally, Theorem \ref{th6.1} reveals that the asymptotic spreading speed (if exists) cannot be faster than the minimal speed for the traveling wave fronts corresponding to the model (\ref{6.1}). It would be very interesting if one can realize how the asymptotic spreading speed depends on these parameters. In \cite{LinG} and \cite{Pan}, Lin and Pan, respectively, have obtained some interesting results for the asymptotic spreading speeds of the model (\ref{6.1}) by constructing appropriate and elaborate upper and/or lower solutions. Their conclusion seems to show that the prey and predator may have different asymptotic spreading speeds.

A great deal of previous mathematical investigation on the spreading of population has been based on the traveling wave fronts of pery-predator system (\ref{6.1}). A striking difference between our present problems and (\ref{6.1}) is that the spreading front in our present problems is given explicitly by a function $x=h(t)$,  beyond which the population densities of both prey and predator are $0$, while in (\ref{6.1}), the two species become positive for all $x$ once $t$ is positive. Secondly, (\ref{6.1}) guarantees successful spreading of the two species for any nontrivial initial populations $(u(0, x)$ and $v(0,x)$, regardless of their initial sizes and supporting area, but the dynamics of our present problems exhibit the spreading-vanishing dichotomy. The phenomenon exhibited by this dichotomy seems closer to the reality.

\vskip 4pt \noindent {\bf Acknowledgment:} The author would like to thank Professor
Guo Lin. He provide me useful references \cite{LinG, LLM, Pan} on the traveling wave fronts and minimal speed of the system (\ref{6.1}). The author is also grateful to the referee for
helpful comments.

\end{document}